\documentclass{basi}
\usepackage{amsmath}
\usepackage{amssymb}
\newcommand{\abc}{{\Bigl(1+\frac{5}{2}A_2\Bigr) }}
\newcommand{\aac}{{(1-\frac{A_2}{2}) }}
\newcommand{\amc}{{(1-\mu) }}
\newcommand{\adc}{{\frac{\delta ^2}{2} }}

\newcommand{\zx}{{(x+\mu)}}
\newcommand{\zox}{{(x+\mu-1)}}
\newcommand{\zd}{{\displaystyle}}

\newcommand{\az}{{1,0}}
\newcommand{\za}{{0,1}}
\newcommand{\zae}{A_2\epsilon}
\newcommand{\zwe}{nW_1\epsilon}
\begin{document}
\title[Higher  Order Normalizations --- Poynting-Robertson
Drag]{Higher  Order Normalizations in the Generalized
Photogravitational Restricted Three Body Problem with
Poynting-Robertson Drag}
\author[B.S.Kushvah, J.P. Sharma and  B.Ishwar ]
       {B.S.Kushvah\thanks{JRF  DST Project, Email:bskush@hotmail.com}, J.P. Sharma\thanks{Co-P.I.   DST Project} and
       B.Ishwar\thanks{P.I. DST Project, Email:ishwar\_ bhola@hotmail.com}\\University Department of Mathematics,\\
    B.R.A. Bihar University Muzaffarpur-842001}
\maketitle
\label{firstpage}
\begin{abstract}
Higher order normalizations are    performed   in the generalized
photogravitational restricted three body problem with
Poynting-Robertson drag. In this problem we have taken bigger
primary as a source of radiation and smaller primary as an oblate
spheroid. Whittaker method is used to transform the second order
part of the Hamiltonian into the normal form.  We have also
performed Birkhoff's normalization of the Hamiltonian. For this we
have utilized Henrard's  method and expanded the coordinates of the
infinitesimal body in double D'Alembert series. We have found the
values of first and second order components. They are affected by
radiation pressure, oblateness and P-R drag. Finally we obtained the
third order part of the Hamiltonian  zero.
\end{abstract}

\section{Introduction}
\label{sec:intro} The restricted three body problem describes the
motion of an infinitesimal mass moving under the gravitational
effect of the two finite masses, called primaries, which move in
circular orbits around their center of mass on account of their
mutual attraction and the infinitesimal mass not influencing the
motion of the primaries. The classical restricted three body problem
is generalized to include the force of radiation pressure, the
Poynting-Robertson effect and oblateness effect.

J. H. Poynting(1903) considered the effect of the absorption and
subsequent re-emission of sunlight by small isolated particles in
the solar system. His work was later modified by H. P. Robertson
(1937) who used  precise relativistic treatments of the first order
in the ratio of the velocity of the particle to that of light.
Chernikov Yu. A.(1970) and Schuerman(1980) who discussed the
position as well as the stability of the Lagrangian equilibrium
points when radiation pressure, P-R drag force are included. Murray
C. D. (1994) systematically discussed the dynamical effect of
general drag in the planar circular restricted three body problem,
Liou J.C.{\it et al.}(1995) examined the effect of radiation
pressure, P-R drag and solar wind drag in the restricted three body
problem.

Moser's conditions(1962), Arnold's theorem(1961) and Liapunov's
theorem (1956) played a significant role in deciding the nonlinear
stability of an equilibrium point. Moser gave some modifications in
Arnold's theorem. Then Deprit and Deprit(1967) investigated the
nonlinear stability of triangular points by applying Moser's
modified version of Arnold's theorem(1961). Maciejewski and
Gozdziewski(1991) described the normalization algorithms of
Hamiltonian near an equilibrium point. Niedzielska(1994)
investigated the nonlinear stability of the libration points in the
photogravitational restricted three body problem. Mishra P. and
Ishwar B.(1995) studied second order normalization in the
generalized restricted problem of three bodies, smaller primary
being an oblate spheroid. Ishwar B.(1997) studied nonlinear
stability in the generalized restricted three body problem.

In this paper higher order normalizations are    performed   in the
generalized photogravitational restricted three body problem with
Poynting-Robertson drag.  Whittaker method is used to transform the
second order part of the Hamiltonian into the normal form. We
 have performed Birkhoff's normalization of the Hamiltonian. For this we
have utilized Henrard's  method and expanded the coordinates of the
third body in double D'Alembert series. We have found the values of
first and second order components. The second order components are
obtained as solutions of the two partial differential equations. We
have employed the first condition of KAM theorem in solving these
equations. The first and second order components are affected by
radiation pressure, oblateness and P-R drag. Finally we obtained the
third order part $H_3$ of the Hamiltonian in $I_1^{1/2}, I_2^{1/2}$
zero.
\section{Location of Triangular Equilibrium Points}
\label{sec:using} Equations of motion are
\begin{align}
\ddot{x}-2n\dot{y}&=U_x ,\quad\text{where},\quad U_x=\frac{\partial{U_1}}{\partial{x}}-\frac{W_{1}N_1}{r^2_1}\\
\ddot{y}+2n\dot{x}&=U_y,\hspace{.85in}U_y=\frac{\partial{U_1}}{\partial{y}}-\frac{W_{1}N_2}{r^2_1}\\
U_1&=\zd{\frac{n^2(x^2+y^2)}{2}}+\frac{\amc{q_1}}{r_1}+\frac{\mu}{r_2}+\frac{\mu{A_2}}{2r^3_2}
\end{align}
\begin{gather*}
r^2_1=\zx^2+y^2,\quad  r^2_2=\zox^2+y^2,\quad n^2=1+\frac{3}{2}A_2,\\
N_1=\frac{\zx[\zx\dot{x}+y\dot{y}]}{r^2_1}+\dot{x}-ny,\quad
N_2=\frac{y[\zx\dot{x}+y\dot{y}]}{r^2_1}+\dot{y}+n\zx
\end{gather*}
$W_1=\frac{(1-\mu)(1-q_1)}{c_d}$,
$\mu=\frac{m_2}{m_1+m_2}\leq\frac{1}{2}$, $m_1,m_2$ be  the  masses
of the primaries, $A_2=\frac{r^2_e-r^2_p}{5r^2}$ be the oblateness
coefficient, $r_e$ and $r_p$ be the equatorial and polar radii
respectively, $r$ be the distance between primaries,
$q=\bigl(1-\frac{F_p}{F_g}\bigr)$ be the mass reduction factor
expressed in terms of the particle's radius $a$, density $\rho$ and
radiation pressure efficiency factor $\chi$ (in the C.G.S.system)
i.e., $q=1-\zd{\frac{5.6\times{10^{-5}}\chi}{a\rho}}$. Assumption
$q=constant$ is equivalent to neglecting fluctuation in the beam of
solar radiation and the effect of solar radiation, the effect of the
planet's shadow, obviously $q\leq1$. Triangular equilibrium points
are given by $U_x=0,U_y=0,z=0,y\neq{0}$, then we have
\begin{align}
 x_*&=x_0\Biggl\{1-\zd{\frac{nW_1\bigl[\amc\abc+\mu\aac\adc\bigr]}{3\mu\amc{y_0 x_0}}}-\adc\frac{A_2}{x_0}\Biggr\} \label{eq:1x}\\
y_*&=y_0{\Biggl\{1-\zd{\frac{nW_1\delta^2\bigl[2\mu-1-\mu(1-\frac{3A_2}{2})\adc+7\amc\frac{A_2}{2}\bigr]}{3\mu\amc{y^3_0}}}-\zd{\frac{\delta^2\bigl(1-\adc)A_2}{y^2_0}}\Biggr\}^{1/2}
}\end{align} where $x_0=\adc-\mu$,
$y_0=\pm\delta\bigl(1-\frac{\delta^2}{4}\bigr)^{1/2}$ and
$\delta=q^{1/3}_1$, as in Kushvah \& Ishwar(2006)
\section{Normalization of $H_2$}
We used  Whittaker (1965) method  for the transformation of $H_2$
into normal form.   The Lagrangian function of the problem can be
written as
\begin{align}
L&=\frac{1}{2}(\dot{x}^2+\dot{y}^2)+n(x\dot{y}-\dot{x}y)+\frac{n^2}{2}(x^2+y^2)+\frac{\amc{q_1}}{r_1}+\frac{\mu}{r_2}+\frac{\mu{A_2}}{2r^3_2}\\\notag
&+W_1\Bigl\{\frac{\zx\dot{x}+y\dot{y}}{2r^2_1}-n
\arctan{\frac{y}{\zx}}\Bigr\}\\\notag
\end{align}
and the Hamiltonian  is $H=-L+p_x\dot{x}+p_y\dot{y}$, where
$p_x,p_y$ are the momenta coordinates given by \[
p_x=\frac{\partial{L}}{\partial{\dot{x}}}=\dot{x}-ny+\frac{W_1}{2r_1^2}\zx,
\quad
p_y=\frac{\partial{L}}{\partial{\dot{y}}}=\dot{y}+nx+\frac{W_1}{2r_1^2}y
\]
For simplicity we suppose  $q_1=1-\epsilon$, with $|\epsilon|<<1$
then coordinates of triangular equilibrium points  can be written in
the form
\begin{align}
x&=\frac{\gamma}{2}-\frac{\epsilon}{3}-\frac{A_2}{2}+\frac{A_2
\epsilon}{3}-\frac{(9+\gamma)}{6\sqrt{3}}nW_1-\frac{4\gamma
\epsilon}{27\sqrt{3}}nW_1 \\
y&=\frac{\sqrt{3}}{2}\Bigl\{1-\frac{2\epsilon}{9}-\frac{A_2}{3}-\frac{2A_2
\epsilon}{9}+\frac{(1+\gamma)}{9\sqrt{3}}nW_1-\frac{4\gamma
\epsilon}{27\sqrt{3}}nW_1\Bigr\}
\end{align}
where $\gamma=1-2\mu$.
 We shift the origin to $L_4$. For that, we change
$x\rightarrow {x_*}+x$ and  $y\rightarrow{y_*}+y$. Let $a=x_*+\mu,
b=y_*$ so that
\begin{align}
a&= \frac{1}{2} \biggl\{1-\frac{2\epsilon}{3}-A_2+\frac{2A_2
\epsilon}{3}-\frac{(9+\gamma)}{3\sqrt{3}}nW_1-\frac{8\gamma
\epsilon}{27\sqrt{3}}nW_1 \biggr\}\\
b&=\frac{\sqrt{3}}{2}\biggl\{1-\frac{2\epsilon}{9}-\frac{A_2}{3}-\frac{2A_2
\epsilon}{9}+\frac{(1+\gamma)}{9\sqrt{3}}nW_1-\frac{4\gamma
\epsilon}{27\sqrt{3}}nW_1\biggr\}
\end{align}
Expanding $L$ in power series of $x $ and $y$, we get
\begin{eqnarray}
 L&=&L_0+L_1+L_2+L_3+\cdots \\
H&=&H_0+H_1+H_2+H_3+\cdots =-L+p_x{\dot{x}}+p_y{\dot{y}}
  \end{eqnarray}
  where $L_0,L_1,L_2,L_3 \ldots$ are
\begin{eqnarray}
L_0&=&\frac{3}{2}-\frac{2\epsilon}{3}-\frac{\gamma
\epsilon}{3}+\frac{ 3 \gamma A_2}{4}-\frac{3 A_2 \epsilon}{2}-\gamma
A_2 \notag \\
&&-\frac{\sqrt{3}nW_1}{4}+\frac{2\gamma nW_1}{3\sqrt{3}}-\frac{ n
\epsilon W_1}{3\sqrt{3}}-\frac{23\epsilon\gamma n W_1}{54\sqrt{3}}-n
W_1\arctan{\frac{b}{a}}
\end{eqnarray}
\begin{align}
L_1&=\dot{x}\biggl\{-\frac{\sqrt{3}}{2}+\frac{\epsilon}{3\sqrt{3}}-\frac{5
A_2 }{8\sqrt{3}}+\frac{7\epsilon A_2}{12\sqrt{3}}+\frac{4
nW_1}{9}-\frac{\gamma nW_1}{18}\biggr\}\notag \\&+
\dot{y}\biggl\{\frac{1}{2}-\frac{\epsilon}{3}-\frac{A_2
}{8}+\frac{\epsilon A_2}{12}-\frac{ \gamma nW_1}{6\sqrt{3}}+\frac{2
\epsilon nW_1}{3\sqrt{3}}\biggr\} \notag\\
& -x \biggl\{-\frac{1}{2}+\frac{\gamma}{2}+\frac{9
A_2}{8}+\frac{15\gamma A_2}{8}-\frac{35\epsilon
A_2}{12}-\frac{29\gamma \epsilon  A_2}{12}+
\frac{3\sqrt{3}nW_1}{8}-\frac{5\epsilon n W_1}{12\sqrt{3}}-\frac{7
\gamma \epsilon nW_1}{4\sqrt{3}}\biggr\}\notag\\&-y
\biggl\{\frac{15\sqrt{3}A_2}{2}+\frac{9\sqrt{3}\gamma
A_2}{8}-2\sqrt{3} \epsilon A_2-2\sqrt{3}\gamma \epsilon A_2-
\frac{nW_1}{8}+\gamma nW_1-\frac{43 \epsilon nW_1}{36}\notag
\\&-\frac{23\gamma\epsilon nW_1}{36}\biggr\}\end{align}
\begin{eqnarray}
L_2&=&\frac{(\dot x^2+ \dot y^2)}{2}+n(x\dot y-\dot x y)+
\frac{n^2}{2}(x^2+y^2)-Ex^2-Fy^2-G xy
\end{eqnarray}
\begin{equation}
L_3=-\frac{1}{3!}\left\{x^3T_1+3x^2yT_2+3xy^2T_3+y^3T_4+6T_5\right\}
  \end{equation}
  where
\begin{eqnarray}
E&=&\frac{1}{16}\Bigl\{ 2-6\epsilon- 3A_2-
\frac{31A_2\epsilon}{2}-\frac{(69+\gamma)}{6\sqrt{3}}nW_1+\frac{2(307+75\gamma)
\epsilon}{27\sqrt{3}}nW_1 \notag \\&+&\gamma \bigl\{2\epsilon+12A_2+
\frac{A_2\epsilon}{3}+\frac{(199+17\gamma)}{6\sqrt{3}}nW_1-\frac{2(226+99\gamma)
\epsilon}{27\sqrt{3}}nW_1\bigr\}\Bigr\}
  \end{eqnarray}
  \begin{eqnarray}
F&=&\frac{-1}{16}\Bigl\{ 10-2\epsilon+21A_2-
\frac{717A_2\epsilon}{18}-\frac{(67+19\gamma)}{6\sqrt{3}}nW_1+\frac{2(413-3\gamma)
\epsilon}{27\sqrt{3}}nW_1 \notag \\&+&\gamma \bigl\{6\epsilon-
\frac{293A_2\epsilon}{18}+\frac{(187+27\gamma)}{6\sqrt{3}}nW_1-\frac{4(247+3\gamma)
\epsilon}{27\sqrt{3}}nW_1\bigr\}\Bigr\}
  \end{eqnarray}
   \begin{eqnarray}
G&=&\frac{\sqrt{3}}{8}\Bigl\{ 2\epsilon+6A_2-
\frac{37A_2\epsilon}{2}-\frac{(13+\gamma)}{2\sqrt{3}}nW_1+\frac{2(79-7\gamma)
\epsilon}{27\sqrt{3}}nW_1 \notag \\&-&\gamma
\bigl\{6-\frac{\epsilon}{3}+13A_2-
\frac{33A_2\epsilon}{2}+\frac{(11-\gamma)}{2\sqrt{3}}nW_1-\frac{(186-\gamma)
\epsilon}{9\sqrt{3}}nW_1\bigr\}\Bigr\}
  \end{eqnarray}
 \begin{eqnarray} T_1&=&\frac{3}{16}\biggl[\frac{16}{3}\epsilon+6A_2-\frac{979}{18}\zae+\frac{(143+9\gamma)}{6\sqrt{3}}nW_1+\frac{(459+376\gamma)}{27\sqrt{3}}\zwe\notag\\&&+\gamma\Biggl\{14+\frac{4\epsilon}{3}+25A_2-\frac{1507 }{18}\zae-\frac{(215+29\gamma)}{6\sqrt{3}}nW_1
\notag
\\&-&\frac{2(1174+169\gamma)}{27\sqrt{3}}\zwe\Biggr\}\biggr]\end{eqnarray}
\begin{eqnarray}
T_2&=&\frac{3\sqrt{3}}{16}\biggl[14-\frac{16}{3}\epsilon+\frac{A_2}{3}-\frac{367}{18}\zae+\frac{115(1+\gamma)}{18\sqrt{3}}nW_1-\frac{(959-136\gamma)}{27\sqrt{3}}\zwe\notag\\&&+\gamma\Biggl\{\frac{32\epsilon}{3}+40A_2-\frac{382}{9}\zae+\frac{(511+53\gamma)}{6\sqrt{3}}nW_1\notag
\\&-&\frac{(2519-24\gamma)}{27\sqrt{3}}\zwe\Biggr\}\biggr]
\end{eqnarray}
\begin{eqnarray} T_3&=&\frac{-9}{16}\biggl[\frac{8}{3}\epsilon+\frac{203A_2}{6}-\frac{625}{54}\zae-\frac{(105+15\gamma)}{18\sqrt{3}}nW_1-\frac{(403-114\gamma)}{81\sqrt{3}}\zwe\notag\\&&+\gamma\Biggl\{2-\frac{4\epsilon}{9}+\frac{55A_2}{2}-\frac{797}{54}\zae+\frac{(197+23\gamma)}{18\sqrt{3}}nW_1\notag \\&-&\frac{(211-32\gamma)}{81\sqrt{3}}\zwe\Biggr\}\biggr]
\end{eqnarray}
\begin{eqnarray} T_4&=&\frac{-9\sqrt{3}}{16}\biggl[2-\frac{8}{3}\epsilon+\frac{23A_2}{3}-44\zae-\frac{(37+\gamma)}{18\sqrt{3}}nW_1-\frac{(219+253\gamma)}{81\sqrt{3}}\zwe\notag\\&&+\gamma\Biggl\{4\epsilon+\frac{88}{27}\zae+\frac{(241+45\gamma)}{18\sqrt{3}}nW_1\notag \\&-&\frac{(1558-126\gamma)}{81\sqrt{3}}\zwe\Biggr\}\biggr]
 \end{eqnarray}
\begin{eqnarray} T_5&=&\frac{W_1}{2(a^2+b^2)^3}\biggl[(a\dot{x}+b\dot{y})\left\{3(ax+by)-(bx-ay)^2\right\}\notag \\&-&2(x\dot{x}+y\dot{y})(ax+by)(a^2+b^2)\biggr]\label{eq:t5}\end{eqnarray}
The second order part $H_2$ of the corresponding Hamiltonian takes
the form
\begin{equation}
H_2=\frac{p_x^2+p_y^2}{2}+n(yp_x-xp_y)+Ex^2+Fy^2+Gxy
\end{equation}
To investigate the stability of the motion, as in Whittaker(1965),
we consider the following set of linear equations in the variables
$x, y$:
 \begin{equation}\label{eq:ax}
 \begin{array}{l c l}
 -\lambda p_x& = & \frac{\partial{H_2}}{\partial x}\\&&\\
 -\lambda p_y& = & \frac{\partial{H_2}}{\partial y}\\
 \text{i.e.}\quad AX&=&0
 \end{array} \quad
 \begin{array}{l c l }
 \lambda x& = & \frac{\partial{H_2}}{\partial p_x}\\&&\\
 \lambda y& = & \frac{\partial{H_2}}{\partial p_y}\\&&
 \end{array}
  \end{equation}
  \begin{equation}
  X=\left[\begin{array}{c}
  x\\
  y\\
  p_x\\
  p_y \end{array}\right] \quad \text{and}
 \quad
 A=\left[\begin{array}{c c c c}
  2E & G&\lambda& -n\\
G&2F&n&\lambda\\
  -\lambda& n& 1& 0\\
  -n & -\lambda& 0& 1\end{array}\right]
  \end{equation}
 Clearly $|A|=0$, implies that the characteristic equation
 corresponding to Hamiltonian $H_2$ is given by
 \begin{equation}
 \lambda^4+2(E+F+n^2)\lambda^2+4EF -G^2+n^4-2n^2(E+F)=0 \label{eq:ch}
 \end{equation}
 This is characteristic equation whose discriminant is
  \begin{equation}
 D=4(E+F+n^2)^2-4\bigl\{4EF-G^2+n^4-2n^2(E+F)\bigr\}
 \end{equation}
 Stability is assured  only when $D>0$.
 i.e
  \begin{eqnarray}
  \mu&<&\mu_{c_0}-0.221895916277307669\epsilon +2.1038871010983331 A_2
  \notag\\&+&
    0.493433373141671349\epsilon A_2 +0.704139054372097028 n W_1 \notag
    \\&+&
    0.401154273957540929 n\epsilon W_1\notag
 \end{eqnarray}
 where $\mu_{c_0}=0.0385208965045513718$.
 When $D>0$ the roots $\pm i \omega_1$ and $\pm i \omega_2$ ($\omega_1,\omega_2$ being the long/short -periodic frequencies) are related to each other  as
 \begin{eqnarray}
   \omega_1^2+\omega_2^2&=& 1-\frac{\gamma \epsilon}{2}+\frac{3\gamma A_2}{2}+\frac{83\epsilon A_2}{12}+\frac{299\gamma\epsilon A_2}{144}-\frac{n W_1}{24\sqrt{3}}+\frac{5 \gamma n W_1}{8\sqrt{3}}-\frac{53 \epsilon n W_1}{54\sqrt{3}}\notag\\
  &&-\frac{5 \gamma^2 n W_1}{24\sqrt{3}}+\frac{173 \gamma \epsilon n W_1}{54\sqrt{3}}-\frac{3 \gamma^2 \epsilon n
  W_1}{36\sqrt{3}}\label{eq:w1+w2} \\
  \omega_1^2\omega_2^2&=&\frac{27}{16} -\frac{27\gamma^2}{16}+\frac{9\epsilon}{8}+\frac{9\gamma\epsilon}{8} -\frac{3\gamma^2\epsilon}{8}+\frac{117\gamma A_2}{16}-\frac{241\epsilon A_2}{32}+\frac{2515\gamma\epsilon A_2}{192}\notag\\\label{eq:w1w2}
  &&+\frac{35n W_1}{16\sqrt{3}}-\frac{55 \sqrt{3}\gamma n W_1}{16}-\frac{5\sqrt{3} \gamma^2 n W_1}{4}-\frac{1277 \epsilon n W_1}{288\sqrt{3}}\notag \\&+&\frac{5021 \gamma \epsilon n W_1}{288\sqrt{3}}+\frac{991 \gamma^2 \epsilon n W_1}{48\sqrt{3}} \\
&&(0<\omega_2<\frac{1}{\sqrt{2}}<\omega_1<1)\notag\end{eqnarray}
From (~\ref{eq:w1+w2}) and  (~\ref{eq:w1w2}) it may be noted that
$\omega_j (j=1,2)$ satisfy
\begin{eqnarray}
  \gamma^2&=& 1+\frac{4\epsilon}{9}-\frac{107\epsilon A_2}{27}+\frac{2\gamma \epsilon }{3}+\frac{1579\gamma\epsilon A_2}{324}-\frac{25nW_1}{27\sqrt{3}}-\frac{55\gamma nW_1}{9\sqrt{3}}+\frac{3809\epsilon nW_1}{486\sqrt{3}}\notag \\&+&\frac{4961\gamma\epsilon nW_1}{486\sqrt{3}}+\biggl(-\frac{16}{27}+\frac{32\epsilon}{243}+\frac{8\gamma\epsilon}{27}+\frac{208 A_2}{81}-\frac{8\gamma A_2}{27}-\frac{4868\epsilon A_2}{729}-\frac{563\gamma\epsilon  A_2}{243}\notag\\&&+\frac{296nW_1}{243\sqrt{3}}-\frac{10\gamma nW_1}{27\sqrt{3}}-\frac{15892\epsilon nW_1}{2187\sqrt{3}}-\frac{1864\gamma\epsilon nW_1}{729\sqrt{3}}\biggr)\omega_j^2\notag\\
  && +\biggl(\frac{16}{27}-\frac{32\epsilon}{243}-\frac{208 A_2}{81}-\frac{1880\epsilon A_2}{729}-\frac{2720nW_1}{2187\sqrt{3}}\notag \\&+&\frac{49552\epsilon nW_1}{6561\sqrt{3}}-\frac{80\gamma\epsilon nW_1}{2187\sqrt{3}}\biggr)\omega_j^4
 \end{eqnarray}
 Alternatively, it can also be seen that if $u=\omega_1\omega_2$,
 then equation (~\ref{eq:w1w2}) gives
\begin{eqnarray}
  \gamma^2&=& 1+\frac{4\epsilon}{9}-\frac{107\epsilon A_2}{27}-\frac{25nW_1}{27\sqrt{3}}+\frac{3809\epsilon nW_1}{486\sqrt{3}}\notag \\&+&\gamma\biggl(\frac{2\epsilon }{3}+\frac{1579\epsilon A_2}{324}-\frac{55\gamma nW_1}{9\sqrt{3}}+\frac{4961\gamma\epsilon
  nW_1}{486\sqrt{3}}\biggr)\notag \\&+&\biggl(-\frac{16}{27}+\frac{32\epsilon}{243}+\frac{208 A_2}{81}-\frac{1880\epsilon A_2}{729}+\frac{320nW_1}{243\sqrt{3}}-\frac{15856\epsilon nW_1}{2187\sqrt{3}}\biggr)u^2 \end{eqnarray}
Following the method for reducing $H_2$ into the normal form, as in
Whittaker(1965), use the transformation
\begin{equation}
 X=JT \quad  \text{where}  \quad X=\left[\begin{array}{c}
x\\y\\p_x\\p_y\end{array}\right],J=[J_{ij}]_{1\leq i,j \leq 4},\
T=\left[\begin{array}{c} Q_1\\Q_2\\P_1\\P_2\end{array}\right]
\end{equation}
\begin{equation}
P_i= (2 I_i\omega_i)^{1/2}\cos{\phi_i},  \quad Q_i= (\frac{2
I_i}{\omega_i})^{1/2}\sin{\phi_i}, \quad (i=1,2)
\end{equation}
The transformation changes the second order part of the Hamiltonian
into the normal form \begin{equation}
H_2=\omega_1I_1-\omega_2I_2\end{equation} The general solution of
the corresponding equations of motion are
\begin{equation}I_i=\text{const.}, \quad \phi_i=\pm \omega_i+\text{const},\  (i=1,2)\label{eq:I}\end{equation}
If the oscillations about $L_4$ are exactly linear, the
Eq.(~\ref{eq:I}) represent the integrals of motion and the
corresponding orbits is  given by
\begin{eqnarray}x&=&J_{13}\sqrt{2\omega_1I_1}\cos{\phi_1}+J_{14}\sqrt{2\omega_2I_2}\cos{\phi_2}\label{eq:xb110}\\
y&=&J_{21}\sqrt{\frac{2I_1}{\omega_1}}\sin{\phi_1}+J_{22}\sqrt{\frac{2I_2}{\omega_2}}\sin{\phi_2}+J_{23}\sqrt{2I_1}{\omega_1}\cos{\phi_1}\notag
\\&+&J_{24}\sqrt{2I_2}{\omega_2}\cos{\phi_2}\label{eq:yb101}\end{eqnarray}
where \begin{eqnarray}
J_{13}&=&\frac{l_1}{2\omega_1k_1}\left\{1-\frac{1}{2l_1^2}\left[\epsilon+\frac{45A_2}{2}-\frac{717A_2\epsilon}{36}+\frac{(67+19\gamma)}{12\sqrt{3}}nW_1
-\frac{(431-3\gamma)}{27\sqrt{3}}nW_1\epsilon\right]\right.\notag\\&&+\frac{\gamma}{2l_1^2}\left[3\epsilon-\frac{29A_2}{36}-\frac{(187+27\gamma)}{12\sqrt{3}}nW_1
-\frac{2(247+3\gamma)}{27\sqrt{3}}nW_1\epsilon\right]\notag\\&&-\frac{1}{2k_1^2}\left[\frac{\epsilon}{2}-3A_2-\frac{73A_2\epsilon}{24}+\frac{(1-9\gamma)}{24\sqrt{3}}nW_1
+\frac{(53-39\gamma)}{54\sqrt{3}}nW_1\epsilon\right]\notag\\&&-\frac{\gamma}{4k_1^2}\left[\epsilon-3A_2-\frac{299A_2\epsilon}{72}-\frac{(6-5\gamma)}{12\sqrt{3}}nW_1
-\frac{(266-93\gamma)}{54\sqrt{3}}nW_1\epsilon\right]\notag\\&&\left.+\frac{\epsilon}{4l_1^2k_1^2}\left[\frac{3A_2}{4}
+\frac{(33+14\gamma)}{12\sqrt{3}}nW_1\right]+\frac{\gamma\epsilon}{8l_1^2k_1^2}\left[\frac{347A_2}{36}
-\frac{(43-8\gamma)}{4\sqrt{3}}nW_1 \right]\right\}\end{eqnarray}
 \begin{eqnarray} J_{14}&=&\frac{l_2}{2\omega_2k_2}\left\{1-\frac{1}{2l_2^2}\left[\epsilon+\frac{45A_2}{2}-\frac{717A_2\epsilon}{36}+\frac{(67+19\gamma)}{12\sqrt{3}}nW_1
-\frac{(431-3\gamma)}{27\sqrt{3}}nW_1\epsilon\right]\right.\notag\\&&-\frac{\gamma}{2l_2^2}\left[3\epsilon-\frac{293A_2}{36}+\frac{(187+27\gamma)}{12\sqrt{3}}nW_1
-\frac{2(247+3\gamma)}{27\sqrt{3}}nW_1\epsilon\right]\notag\\&&-\frac{1}{2k_2^2}\left[\frac{\epsilon}{2}-3A_2-\frac{73A_2\epsilon}{24}+\frac{(1-9\gamma)}{24\sqrt{3}}nW_1
+\frac{(53-39\gamma)}{54\sqrt{3}}nW_1\epsilon\right]\notag\\&&+\frac{\gamma}{2k_2^2}\left[\epsilon-3A_2-\frac{299A_2\epsilon}{72}-\frac{(6-5\gamma)}{12\sqrt{3}}nW_1
-\frac{(268-9\gamma)}{54\sqrt{3}}nW_1\epsilon\right]\notag\\&&\left.-\frac{\epsilon}{4l_2^2k_2^2}\left[\frac{33A_2}{4}
+\frac{(1643-93\gamma)}{216\sqrt{3}}nW_1\right]+\frac{\gamma\epsilon}{4l_2^2k_2^2}\left[\frac{737A_2}{72}
-\frac{(13+2\gamma)}{\sqrt{3}}nW_1 \right]\right\}\end{eqnarray}
\begin{eqnarray} J_{21}&=&-\frac{4n\omega_1}{l_1k_1}\left\{1+\frac{1}{2l_1^2}\left[\epsilon+\frac{45A_2}{2}-\frac{717A_2\epsilon}{36}+\frac{(67+19\gamma)}{12\sqrt{3}}nW_1
-\frac{(413-3\gamma)}{27\sqrt{3}}nW_1\epsilon\right]\right.\notag\\&&-\frac{\gamma}{2l_1^2}\left[3\epsilon-\frac{293A_2}{36}+\frac{(187+27\gamma)}{12\sqrt{3}}nW_1
-\frac{2(247+3\gamma)}{27\sqrt{3}}nW_1\epsilon\right]\notag\\&&-\frac{1}{2k_1^2}\left[\frac{\epsilon}{2}-3A_2-\frac{73A_2\epsilon}{24}+\frac{(1-9\gamma)}{24\sqrt{3}}nW_1
+\frac{(53-39\gamma)}{54\sqrt{3}}nW_1\epsilon\right]\notag\\&&-\frac{\gamma}{4k_1^2}\left[\epsilon-3A_2-\frac{299A_2\epsilon}{72}-\frac{(6-5\gamma)}{12\sqrt{3}}nW_1
-\frac{(268-93\gamma)}{54\sqrt{3}}nW_1\epsilon\right]\notag\\&&\left.+\frac{\epsilon}{8l_1^2k_1^2}\left[\frac{33A_2}{4}+\frac{(68-10\gamma)}{24\sqrt{3}}nW_1\right]+\frac{\gamma\epsilon}{8l_1^2k_1^2}\left[\frac{242A_2}{9}
+\frac{(43-8\gamma)}{4\sqrt{3}}nW_1 \right]\right\}\end{eqnarray}
\begin{eqnarray} J_{22}&=&\frac{4n\omega_2}{l_2k_2}\left\{1+\frac{1}{2l_2^2}\left[\epsilon+\frac{45A_2}{2}-\frac{717A_2\epsilon}{36}+\frac{(67+19\gamma)}{12\sqrt{3}}nW_1
-\frac{(413-3\gamma)}{27\sqrt{3}}nW_1\epsilon\right]\right.\notag\\&&-\frac{\gamma}{2l_2^2}\left[3\epsilon-\frac{293A_2}{36}+\frac{(187+27\gamma)}{12\sqrt{3}}nW_1
-\frac{2(247+3\gamma)}{27\sqrt{3}}nW_1\epsilon\right]\notag\\&&+\frac{1}{2k_2^2}\left[\frac{\epsilon}{2}-3A_2-\frac{73A_2\epsilon}{24}+\frac{(1-9\gamma)}{24\sqrt{3}}nW_1
+\frac{(53-39\gamma)}{54\sqrt{3}}nW_1\epsilon\right]\notag\\&&-\frac{\gamma}{4k_2^2}\left[\epsilon-3A_2-\frac{299A_2\epsilon}{72}-\frac{(6-5\gamma)}{12\sqrt{3}}nW_1
-\frac{(268-93\gamma)}{54\sqrt{3}}nW_1\epsilon\right]\notag\\&&\left.+\frac{\epsilon}{4l_2^2k_2^2}\left[\frac{33A_2}{4}+\frac{(34+5\gamma)}{12\sqrt{3}}nW_1
\right]+\frac{\gamma\epsilon}{8l_2^2k_2^2}\left[\frac{75A_2}{2}+\frac{(43-8\gamma)}{4\sqrt{3}}nW_1
 \right]\right\}\end{eqnarray}
\begin{eqnarray} J_{23}&=&\frac{\sqrt{3}}{4\omega_1l_1k_1}\left\{2\epsilon+6A_2+\frac{37A_2\epsilon}{2}-\frac{(13+\gamma)}{2\sqrt{3}}nW_1
+\frac{2(79-7\gamma)}{9\sqrt{3}}nW_1\epsilon\right.\notag\\&&-\gamma\left[6+\frac{2\epsilon}{3}+13A_2-\frac{33A_2\epsilon}{2}+\frac{(11-\gamma)}{2\sqrt{3}}nW_1
-\frac{(186-\gamma)}{9\sqrt{3}}nW_1\epsilon\right]\notag\\&&+\frac{1}{2l_1^2}\left[51A_2+\frac{(14+8\gamma)}{3\sqrt{3}}nW_1\right]-\frac{\epsilon}{k_1^2}\left[3A_2
+\frac{(19+6\gamma)}{6\sqrt{3}}nW_1\right]\notag\\&&-\frac{\gamma}{2l_1^2}\left[6\epsilon+135A_2-\frac{808A_2\epsilon}{9}-\frac{(67+19\gamma)}{2\sqrt{3}}nW_1
-\frac{(755+19\gamma)}{9\sqrt{3}}nW_1\epsilon\right]\notag\\&&-\frac{\gamma}{2k_1^2}\left[3\epsilon-18A_2-\frac{55A_2\epsilon}{4}-\frac{(1-9\gamma)}{4\sqrt{3}}nW_1
+\frac{(923-60\gamma)}{12\sqrt{3}}nW_1\epsilon\right]\notag\\&&\left.+\frac{\gamma\epsilon}{8l_1^2k_1^2}\left[\frac{9A_2}{2}
+\frac{(34-5\gamma)}{2\sqrt{3}}nW_1\right]\right\}\qquad\end{eqnarray}
\begin{eqnarray}                 J_{24}&=&\frac{\sqrt{3}}{4\omega_2l_2k_2}\left\{2\epsilon+6A_2+\frac{37A_2\epsilon}{2}-\frac{(13+\gamma)}{2\sqrt{3}}nW_1
+\frac{2(79-7\gamma)}{9\sqrt{3}}nW_1\epsilon\right.\notag\\
&&-\gamma\left[6+\frac{2\epsilon}{3}+13A_2-\frac{33A_2\epsilon}{2}+\frac{(11-\gamma)}{2\sqrt{3}}nW_1
-\frac{(186-\gamma)}{9\sqrt{3}}nW_1\epsilon\right]\notag\\&&-\frac{1}{2l_2^2}\left[51A_2+\frac{(14+8\gamma)}{3\sqrt{3}}nW_1\right]-\frac{\epsilon}{k_2^2}\left[3A_2
+\frac{(19+6\gamma)}{6\sqrt{3}}nW_1\right]\notag\\&&-\frac{\gamma}{2l_2^2}\left[6\epsilon+135A_2-\frac{808A_2\epsilon}{9}-\frac{(67+19\gamma)}{2\sqrt{3}}nW_1
-\frac{(755+19\gamma)}{9\sqrt{3}}nW_1\epsilon\right]\notag\\&&-\frac{\gamma}{2k_1^2}\left[3\epsilon-18A_2-\frac{55A_2\epsilon}{4}-\frac{(1-9\gamma)}{4\sqrt{3}}nW_1
+\frac{(923-60\gamma)}{12\sqrt{3}}nW_1\epsilon\right]\notag\\&&\left.-\frac{\gamma\epsilon}{4l_1^2k_1^2}\left[\frac{99A_2}{2}
+\frac{(34-5\gamma)}{2\sqrt{3}}nW_1\right]\right\}\label{eq:j24}\end{eqnarray}
with $l_j^2=4\omega_j^2+9,(j=1,2)$ and $ k_1^2=2\omega_1^2-1,
k_2^2=-2\omega_2^2+1 $.
\section{Second Order Normalization}
In order to perform Birkhoff's normalization, we use Henrard's
  method[Deprit and Deprit Bartholom\'{e} (1967)] for which the
  coordinates $(x,y)$ of infinitesimal body, to be expanded in
  double D'Alembert series\   $x=\sum_{n\geq1}B_n^\az,$ $y=\sum_{n\geq 1}B_n^\za$
  where the homogeneous components $B_n^\az$ and $B_n^\za$ of
  degree $n$ are of the form
  \begin{equation}
\sum_{0\leq{m}\leq{n}}
I_1^{\frac{n-m}{2}}I_2^{\frac{m}{2}}\sum_{(p,q)}C_{n-m,m,p,q}
\cos{(p\phi_1+q\phi_2)}+S_{n-m,m,p,q} \sin{(p\phi_1+q\phi_2)}
  \end{equation}
  The conditions in double summation are (i) $p$ runs over those
  integers in the interval $0\leq p\leq n-m$ that have the same
  parity as $n-m$ (ii) $q$ runs over those integers in the interval $-m\leq q\leq
  m$ that have the same parity as $m$. Here $I_1$, $I_2$ are the
  action momenta coordinates which are to be taken as constants of
  integer, $\phi_1$, $\phi_2$ are angle coordinates to be
  determined as linear functions of time in such a way that $\dot\phi_1=\omega_1+\sum_{n\geq 1}f_{2n}(I_1,I_2),\dot\phi_2=-\omega_2+\sum_{n\geq 1}g_{2n}(I_1,I_2)$  where  $\omega_1,\omega_2$ are the basic  frequencies, $f_{2n}$ and  $g_{2n}$ are of the form
  \begin{eqnarray}
  f_{2n}&=&\sum_{0\leq m\leq n}{f'}_{2(n-m),2m}I_1^{n-m}I_2^m\\
  g_{2n}&=&\sum_{0\leq m\leq n}{g'}_{2(n-m),2m}I_1^{n-m}I_2^m
  \end{eqnarray}
The first order components $B_1^\az$ and $B_1^\za$ are the values of
 $x$ and  $y$ given by (~\ref{eq:xb110}) (~\ref{eq:yb101}).
In order to findout  the second order components $B_2^\az,B_2^\za$
we consider  Lagrange's  equations of motion
\begin{equation}
\frac{d}{dt}(\frac{\partial L}{\partial \dot x })-\frac{\partial
L}{\partial x }=0, \quad \frac{d}{dt}(\frac{\partial L}{\partial
\dot y })-\frac{\partial L}{\partial y }=0 \end{equation}
\begin{equation}
\text{i.e.}\quad \left.\begin{array}{l c l}
\ddot x-2n\dot y+(2E-n^2)x+Gy&=&\frac{\partial L_3}{\partial x }+\frac{\partial L_4}{\partial x }\\
&&\\\ddot y+2n\dot x+(2F-n^2)y+Gx&=&\frac{\partial L_3}{\partial y
}+\frac{\partial L_4}{\partial y }\end{array}
\right\}\label{eq:lgeq}\end{equation} Since $x$ and $y$ are double
D'Alembert series, $x^jx^k(j\geq0,k\geq0,j+k\geq0)$ and  the time
derivatives $\dot x ,\dot y ,\ddot x, \ddot y $ are also double
D'Alembert series. We can write
\[\dot x=\sum_{n\geq 1} \dot x_n, \quad\dot y=\sum_{n\geq 1} \dot
y_n,\quad\ddot x=\sum_{n\geq 1} \ddot x_n,\quad \ddot y=\sum_{n\geq
1} \ddot y_n \] where $\dot x ,\dot y ,\ddot x, \ddot y $ are
homogeneous components of degree $n$ in $I_1^{1/2},I_2^{1/2}$ i.e.
\begin{eqnarray} \dot x &=&
\frac{d}{dt}\sum_{n\geq 1}B_n^\az=\sum_{n\geq
1}\Biggl[\frac{\partial
B_n^\az}{\partial{\phi_1}}(\omega_1+f_2+f_4+\cdots)\notag
\\&+&\frac{\partial
B_n^\az}{\partial{\phi_2}}(-\omega_2+g_2+g_4+\cdots)\Biggr]\end{eqnarray}
We write three components $\dot x_1 ,\dot x_2 ,\dot x_3$ of $\dot x$
\begin{eqnarray}
\dot x_1&=&\omega_1\frac{\partial
B_1^\az}{\partial{\phi_1}}-\omega_2\frac{\partial
B_1^\az}{\partial{\phi_2}}=DB_1^\az\\
 \dot
x_2&=&\omega_1\frac{\partial
B_2^\az}{\partial{\phi_1}}-\omega_2\frac{\partial
B_2^\az}{\partial{\phi_2}}=DB_2^\az\\
\dot x_3&=&\omega_1\frac{\partial
B_3^\az}{\partial{\phi_1}}-\omega_2\frac{\partial
B_3^\az}{\partial{\phi_2}}+f_2\frac{\partial
B_1^\az}{\partial{\phi_1}}-g_2\frac{\partial
B_1^\az}{\partial{\phi_2}}\notag\\
&=&DB_3^\az+f_2\frac{\partial
B_1^\az}{\partial{\phi_1}}-g_2\frac{\partial
B_1^\az}{\partial{\phi_2}}
\end{eqnarray}
where \begin{equation}D\equiv \omega_1\frac{\partial\
}{\partial{\phi_1}}-\omega_2\frac{\partial\
}{\partial{\phi_2}}\end{equation}
 Similarly three components $\ddot
x_1 ,\ddot x_2 ,\ddot x_3$ of $\ddot x$ are
\begin{eqnarray*}
\ddot x_1 &=&D^2B_1^\az, \quad \ddot x_2=D^2B_2^\az,\quad \ddot
x_3=D^2B_3^\az+2\omega_1f_2\frac{\partial^2B_1^\az}{\partial\phi_1^2}-2\omega_2g_2\frac{\partial^2B_1^\az}{\partial\phi_2^2}
\end{eqnarray*}
In similar manner we can write the components of $\dot y, \ddot y$.
Putting the values of  $x, y, \dot x ,\dot y ,\ddot x $ and  $\ddot
y$ in terms of double D'Alembert series in equation (~\ref{eq:lgeq})
we get
\begin{equation}
\left(D^2+2E-1-\frac{3}{2}A_2\right)B_2^\az-\left\{2\left(
1+\frac{3}{4}A_2\right)D-G\right\}B_2^\za=X_2 \label{eq:x2}
\end{equation}
\begin{equation}\label{eq:y2}
\left\{2\left(
1+\frac{3}{4}A_2\right)D+G\right\}B_2^\az+\left(D^2+2F-1-\frac{3}{2}A_2\right)B_2^\za=Y_2
\end{equation} where \[X_2=\left[\frac{\partial
L_3}{\partial x}\right]_{x=B_1^\az,y=B_1^\za} \quad \text{and} \quad
Y_2=\left[\frac{\partial L_3}{\partial
y}\right]_{x=B_1^\az,y=B_1^\za}\] These are two simultaneous partial
differential equations in $B_2^\az$ and $B_2^\za$. We solve these
equations to find the values of $B_2^\az$ and $B_2^\za$, from Eq.
(~\ref{eq:x2}) and (~\ref{eq:y2})
\begin{equation}
\triangle_1 \triangle_2B_2^\az=\Phi_2, \quad \triangle_1
\triangle_2B_2^\za=-\Psi_2 \label{eq:phi_si} \quad \text{where}
\quad \triangle_1=D^2+\omega_1^2, \triangle_2=D^2+\omega_2^2
\end{equation}
\begin{equation}
\Phi_2=(D^2+2F-n^2)X_2+(2nD-G)Y_2 \label{eq:phi2}
\end{equation}
\begin{equation}
\Psi_2=(2nD+G)X_2-(D^2+2E-n^2)Y_2 \label{eq:psi2}
\end{equation}
The Eq.(~\ref{eq:phi_si}) can  be solved for $B_2^\az$ and $B_2^\za$
by putting the formula
\[\frac{1}{\triangle_1\triangle_2}\left\{\begin{array}{c}\cos(p\phi_1+q\phi_2)\\ \mbox{or} \\\sin(p\phi_1+q\phi_2)\end{array}=\frac{1}{\triangle_{p,q}}\left\{\begin{array}{c}\cos(p\phi_1+q\phi_2)\\\mbox{or} \\\sin(p\phi_1+q\phi_2)\end{array}\right.\right.\]
where \[\triangle_{p,q}=\left[
\omega_1^2-(\omega_1p-\omega_2q)^2\right]\left[
\omega_2^2-(\omega_1p-\omega_2q)^2\right]
\]
provided $\triangle_{p,q}\neq0$. Since $\triangle_{1,0}=0,
\triangle_{0,1}=0$ the terms
$\cos\phi_1,\sin\phi_1,\cos\phi_2,\sin\phi_2$ are the critical
terms, $\Phi_2$ and $\Psi_2$ are free from such terms. By
condition(1) of Moser's theorem $k_1\omega_1+k_2\omega_2\neq 0$  for
all pairs $(k_1,k_2)$ of integers such that $|k_1|+|k_2|\leq4$,
therefore each of $\omega_1, \omega_2,
\omega_1\pm2\omega_2,\omega_2\pm2\omega_1$ is different from zero
and consequently none of the divisors $\triangle_{0,0},
\triangle_{0,2}, \triangle_{2,0}, \triangle_{1,1}, \triangle_{1,-1}$
is zero. The second order components $B_2^\az, B_2^\za$ are as
follows:
\begin{eqnarray}
B_2^\az&=&r_1I_1+r_2I_2+r_3I_1\cos2\phi_1+r_4I_2\cos2\phi_2+r_5I_1^{1/2}I_2^{1/2}\cos(\phi_1-\phi_2)\notag\\&&
+r_6I_1^{1/2}I_2^{1/2}\cos(\phi_1+\phi_2)+r_7I_1\sin2\phi_1+r_8I_2\sin2\phi_2\notag\\
&&+r_9I_1^{1/2}I_2^{1/2}\sin(\phi_1-\phi_2)+r_{10}I_1^{1/2}I_2^{1/2}\sin(\phi_1+\phi_2)\label{eq:b2az}
\end{eqnarray} and\begin{eqnarray}
B_2^\za&=&-\left\{s_1I_1+s_2I_2+s_3I_1\cos2\phi_1+s_4I_2\cos2\phi_2+s_5I_1^{1/2}I_2^{1/2}\cos(\phi_1-\phi_2)\right.\notag\\
&&
+s_6I_1^{1/2}I_2^{1/2}\cos(\phi_1+\phi_2)+s_7I_1\sin2\phi_1+s_8I_2\sin2\phi_2\notag\\
&&+\left.s_9I_1^{1/2}I_2^{1/2}\sin(\phi_1-\phi_2)+s_{10}I_1^{1/2}I_2^{1/2}\sin(\phi_1+\phi_2)\right\}\label{eq:b2za}
\end{eqnarray}
where
\begin{eqnarray}
&r_1=&\frac{1}{\omega_1^2\omega_2^2}\left\{
J_{13}^2\omega_1F_4+J_{13}J_{23}\omega_1F_4'+\left(\frac{J_{21}^2}{\omega_1}+J_{23}^2\omega_1\right)F_4''\right\}
\end{eqnarray}
\begin{eqnarray}
&r_2=&\frac{1}{\omega_1^2\omega_2^2}\left\{
J_{14}^2\omega_2F_4+J_{14}J_{24}\omega_2F_4'+\left(\frac{J_{22}^2}{\omega_2}+J_{24}^2\omega_2\right)F_4{''}\right\}
\end{eqnarray}
\begin{eqnarray}
&r_3=&\frac{-1}{3\omega_1^2(4\omega_1^2-\omega_2^2)}\Biggl\{
8\omega_1^3J_{21}(J_{13}F_1'+2J_{23}F_1'')+4\omega_1^2\biggl[(J_{13}F_2+J_{23}F_2'')J_{13}\omega_1\notag
\\&-&\biggl(\frac{J_{21}^2}{\omega_1}-J_{23}^2\omega_1\biggr)F_1''\biggr]-2\omega_1J_{21}(J_{13}F_3'+2J_{23}F_3'')-\omega_1J_{13}(J_{13}F_4+J_{23}F_4'')\omega_1
\notag
\\&+&\biggl(\frac{J_{21}^2}{\omega_1}-J_{23}^2\omega_1\biggr)F_1''\Biggr\}\end{eqnarray}
\begin{eqnarray}
&r_4=&\frac{1}{3\omega_2^2(4\omega_2^2-\omega_1^2)}\Biggl\{8\omega_2^3J_{22}(J_{14}F_1'+2J_{24}F_1'')-4\omega_2^2\biggl[(J_{14}F_2+J_{24}F_2'')J_{14}\omega_2\notag
\\&-&\biggl(\frac{J_{22}^2}{\omega_2}-J_{24}^2\omega_2\biggr)F_2''\biggr]-2\omega_2J_{22}(J_{14}F_3'+2J_{24}F_3'')-\omega_2J_{14}(J_{14}F_4+J_{24}F_4'')\omega_2\notag \\&-&\biggl(\frac{J_{22}^2}{\omega_2}-J_{24}^2\omega_2\biggr)F_4''\Biggr\}\end{eqnarray}
\begin{eqnarray}
&r_5&=\frac{1}{\omega_1\omega_2(2\omega_1+\omega_2)(4\omega_1+2\omega_2)}\Biggl\{(\omega_1+\omega_2)^3\biggl[\bigl\{J_{13}J_{22}(\frac{\omega_1}{\omega_2})^{1/2}-J_{14}J_{21}(\frac{\omega_2}{\omega_1})^{1/2}\bigr\}F_1'\notag\\&&-2\bigl\{J_{21}J_{24}(\frac{\omega_2}{\omega_1})^{1/2}-J_{22}J_{23}(\frac{\omega_1}{\omega_2})^{1/2}\bigr\}F_1''\biggl]-(\omega_1+\omega_2)^2\biggl[\bigr\{2\bigl\{J_{13}J_{14}F_2\notag\\&&+(J_{13}J_{24}+J_{14}J_{23})F_2'\bigr\}(\omega_1\omega_2)^{1/2}+\bigr\{\frac{J_{21}J_{22}}{(\omega_1\omega_2)^{1/2}}+J_{23}J_{24}(\omega_1\omega_2)^{1/2}\bigr\}F_2''\biggr]\notag\\
&&-(\omega_1+\omega_2)\biggl[\bigl\{J_{13}J_{22}(\frac{\omega_1}{\omega_2})^{1/2}-J_{14}J_{21}(\frac{\omega_2}{\omega_1})^{1/2}\bigr\}F_3'-2\bigl\{J_{21}J_{24}(\frac{\omega_2}{\omega_1})^{1/2}\notag
\\&-&J_{22}J_{23}(\frac{\omega_1}{\omega_2})^{1/2}\bigr\}F_3''\biggl]+\biggl[\bigr\{2\bigl\{J_{13}J_{14}F_4+(J_{13}J_{24}+J_{14}J_{23})F_4'\bigr\}(\omega_1\omega_2)^{1/2}\notag \\&+&2\bigr\{\frac{J_{21}J_{22}}{(\omega_1\omega_2)^{1/2}}+J_{23}J_{24}(\omega_1\omega_2)^{1/2}\bigr\}F_4''\biggr]\Biggr\}
\end{eqnarray}
\begin{eqnarray}
&r_6=&\frac{-1}{\omega_1\omega_2(2\omega_1-\omega_2)(4\omega_1-2\omega_2)}
\Biggl\{(\omega_1-\omega_2)^3\biggl[\bigl\{J_{13}J_{22}(\frac{\omega_1}{\omega_2})^{1/2}-J_{14}J_{21}(\frac{\omega_2}{\omega_1})^{1/2}\bigr\}F_1'\notag\\
&&+2\bigl\{J_{21}J_{24}(\frac{\omega_2}{\omega_1})^{1/2}+J_{22}J_{23}(\frac{\omega_1}{\omega_2})^{1/2}\bigr\}F_1''\biggl]\notag\\
&&+(\omega_1-\omega_2)^2\biggl[\bigr\{2\bigl\{J_{13}J_{14}F_2+(J_{13}J_{24}+J_{14}J_{23})F_2'\bigr\}(\omega_1\omega_2)^{1/2}\notag\\
&&-2\bigr\{\frac{J_{21}J_{22}}{(\omega_1\omega_2)^{1/2}}-J_{23}J_{24}(\omega_1\omega_2)^{1/2}\bigr\}F_2''\biggr]-(\omega_1-\omega_2)\biggl[\bigl\{J_{13}J_{22}(\frac{\omega_1}{\omega_2})^{1/2}\notag\\
&&-J_{14}J_{21}(\frac{\omega_2}{\omega_1})^{1/2}\bigr\}F_3'+2\bigl\{J_{21}J_{22}(\frac{\omega_2}{\omega_1})^{1/2}+J_{22}J_{23}(\frac{\omega_1}{\omega_2})^{1/2}\bigr\}F_3''\biggl]\notag\\
&&-\biggl[\bigr\{2\bigl\{J_{13}J_{14}F_4+(J_{13}J_{24}+J_{14}J_{23})F_4'\bigr\}(\omega_1\omega_2)^{1/2}\notag
\\&&-2\bigr\{\frac{J_{21}J_{22}}{(\omega_1\omega_2)^{1/2}}-J_{23}J_{24}(\omega_1\omega_2)^{1/2}\bigr\}F_4''\biggr]\Biggr\}\notag\\&&
\end{eqnarray}
\begin{eqnarray}
&r_7=&\frac{1}{3\omega_1^2(4\omega_1^2-\omega_2^2)} \Biggl\{
8\omega_1^3\biggl[J_{13}(J_{13}F_1+J_{23}F_1')\omega_1-\biggl(\frac{J_{21}^2}{\omega_1}-J_{23}^2\omega_1\biggr)F_1''\biggr]\notag\\&&-2\omega_1\biggl[\omega_1J_{13}(J_{13}F_3+J_{23}F_3')-\biggl(\frac{J_{21}^2}{\omega_1}-J_{23}^2\omega_1\biggr)F_3''\biggr]\notag\\&&-4\omega_1^2J_{21}(J_{13}F_2+J_{23}F_2'')\omega_1+J_{21}(J_{13}F_4'+2J_{23}F_4'')\Biggr\}
\end{eqnarray}
\begin{eqnarray}
&r_8&=\frac{-1}{3\omega_2^2(4\omega_2^2-\omega_1^2)} \Biggl\{
8\omega_2^3\biggl[J_{14}(J_{14}F_1+J_{24}F_1')\omega_2-\biggl(\frac{J_{22}^2}{\omega_2}-J_{24}^2\omega_2\biggr)F_1''\biggr]\notag\\
&&+4\omega_2^2J_{22}(J_{14}F_2+2J_{24}F_2'')\omega_2-2\omega_2\biggl[\omega_2J_{14}(J_{14}F_3+J_{24}F_3')\notag\\
&&-\biggl(\frac{J_{22}^2}{\omega_2}-J_{24}^2\omega_2\biggr)F_3''\biggr]-J_{22}(J_{14}F_4'+2J_{24}F_4'')\Biggr\}
\end{eqnarray}
\begin{eqnarray}
&r_9=&\frac{1}{\omega_1\omega_2(2\omega_1+\omega_2)(\omega_1+2\omega_2)}
\Biggl\{(\omega_1+\omega_2)^3\biggl[\bigr\{2J_{13}J_{14}F_1\notag\\
&&+(J_{13}J_{24}+J_{14}J_{23})F_1'\bigr\}(\omega_1\omega_2)^{1/2}+2\bigr\{\frac{J_{21}J_{22}}{(\omega_1\omega_2)^{1/2}}+J_{23}J_{24}(\omega_1\omega_2)^{1/2}\bigr\}F_1''\biggr]\notag\\
&&-(\omega_1+\omega_2)^2\biggl[\bigl\{J_{13}J_{22}(\frac{\omega_1}{\omega_2})^{1/2}-J_{14}J_{21}(\frac{\omega_2}{\omega_1})^{1/2}\bigr\}F_2'\notag\\
&&-2\bigl\{J_{21}J_{24}(\frac{\omega_2}{\omega_1})^{1/2}-J_{22}J_{23}(\frac{\omega_1}{\omega_2})^{1/2}\bigr\}F_2''\biggl]
\notag\\
&&-(\omega_1+\omega_2)\biggl[\bigr\{2\bigl\{J_{13}J_{14}F_3+(J_{13}J_{24}+J_{14}J_{23})F_3'\bigr\}(\omega_1\omega_2)^{1/2}\notag\\
&&+2\bigr\{\frac{J_{21}J_{22}}{(\omega_1\omega_2)^{1/2}}+J_{23}J_{24}(\omega_1\omega_2)^{1/2}\bigr\}F_3''\biggr]-\biggl[\bigl\{J_{13}J_{22}(\frac{\omega_1}{\omega_2})^{1/2}\notag\\
&&-J_{14}J_{21}(\frac{\omega_2}{\omega_1})^{1/2}\bigr\}F_4'
-2\bigl\{J_{21}J_{24}(\frac{\omega_2}{\omega_1})^{1/2}-J_{22}J_{23}(\frac{\omega_1}{\omega_2})^{1/2}\bigr\}F_4''\biggl]
\Biggr\}
\end{eqnarray}
\begin{eqnarray}
 & r_{10}=&\frac{1}{\omega_1\omega_2(2\omega_1-\omega_2)(2\omega_2-\omega_1)}
\Biggl\{(\omega_1-\omega_2)^3\biggl[\bigr\{2J_{13}J_{14}F_1\notag\\
&&+(J_{13}J_{24}+J_{14}J_{23})F_1'\bigr\}(\omega_1\omega_2)^{1/2}-2\bigr\{\frac{J_{21}J_{22}}{(\omega_1\omega_2)^{1/2}}-J_{23}J_{24}(\omega_1\omega_2)^{1/2}\bigr\}F_1''\biggr]\notag\\
&&-(\omega_1-\omega_2)^2\biggl[\bigl\{J_{13}J_{22}(\frac{\omega_1}{\omega_2})^{1/2}-J_{14}J_{21}(\frac{\omega_2}{\omega_1})^{1/2}\bigr\}F_2'\notag
\\&&+2\bigl\{J_{21}J_{24}(\frac{\omega_2}{\omega_1})^{1/2}+J_{22}J_{23}(\frac{\omega_1}{\omega_2})^{1/2}\bigr\}F_2''\biggl]
\notag\\
&&-(\omega_1-\omega_2)\biggl[\bigr\{2\bigl\{J_{13}J_{14}F_3+(J_{13}J_{24}+J_{14}J_{23})F_3'\bigr\}(\omega_1\omega_2)^{1/2}\notag \\&&-2\bigr\{\frac{J_{21}J_{22}}{(\omega_1\omega_2)^{1/2}}-J_{23}J_{24}(\omega_1\omega_2)^{1/2}\bigr\}F_3''\biggr]\notag\\
&&+\biggl[\bigl\{J_{13}J_{22}(\frac{\omega_1}{\omega_2})^{1/2}-J_{14}J_{21}(\frac{\omega_2}{\omega_1})^{1/2}\bigr\}F_4'\notag
\\&&+2\bigl\{J_{21}J_{24}(\frac{\omega_2}{\omega_1})^{1/2}-J_{22}J_{23}(\frac{\omega_1}{\omega_2})^{1/2}\bigr\}F_4''\biggl]
\Biggr\}
\end{eqnarray}
We can write expressions of $s_i$ with the help of  $r_i$  replacing
$F_i$ by $G_i$, $F_i'$ by $G_i'$  and $F_i''$ by
$G_i''$,$(i=1,2,3,4)$, where
\begin{eqnarray}
 F_1&=&\frac{-\zwe}{6}
\end{eqnarray}
\begin{eqnarray}
  F_2&=&\frac{3}{32}\biggl[\frac{16}{3}\epsilon+6A_2-\frac{979}{18}\zae+\frac{(143+9\gamma)}{6\sqrt{3}}nW_1+\frac{(555+376\gamma)}{27\sqrt{3}}\zwe\notag\\&&+\gamma\Biggl\{14+\frac{4\epsilon}{3}+25A_2-\frac{1507 }{18}\zae-\frac{(215+29\gamma)}{6\sqrt{3}}nW_1
\notag\\
&&-\frac{2(1174+169\gamma)}{27\sqrt{3}}\zwe\Biggr\}\biggr]
\end{eqnarray}
\begin{eqnarray}
F_3&=&\frac{3\sqrt{3}}{16}\biggl[14-\frac{16}{3}\epsilon+\frac{23A_2}{2}-\frac{104}{9}\zae+\frac{115(1+\gamma)}{18\sqrt{3}}nW_1-\frac{2(439-68\gamma)}{27\sqrt{3}}\zwe\notag\\&&+\gamma\Biggl\{\frac{32\epsilon}{3}+40A_2-\frac{310
}{9}\zae+\frac{(511+53\gamma)}{6\sqrt{3}}nW_1
-\frac{(2519-249\gamma)}{27\sqrt{3}}\zwe\Biggr\}\biggr]
\end{eqnarray}
\begin{eqnarray}
F_4&=&\frac{-3}{256}\biggl[364+420A_2-\frac{17801A_2}{9}\zae+\frac{(2821+189\gamma)}{3\sqrt{3}}nW_1-\frac{(23077+9592\gamma)}{27\sqrt{3}}\zwe\notag\\&&+28\gamma\Biggl\{23+\frac{100\epsilon}{21}+\frac{849A_2}{14}+\frac{59
}{7}\zae-\frac{(125+38\gamma)}{6\sqrt{3}}nW_1
\notag\\
&&-\frac{(87613-213\gamma)}{27\sqrt{3}}\zwe\Biggr\}\biggr]
\end{eqnarray}
\begin{eqnarray}
 F_1'&=&\frac{\zwe}{3\sqrt{3}}
\end{eqnarray}
\begin{eqnarray} F_2'&=&\frac{3\sqrt{3}}{16}\biggl[14-\frac{16}{3}\epsilon+A_2-\frac{1367}{18}\zae+\frac{115(1+\gamma)}{18\sqrt{3}}nW_1-\frac{(863-136\gamma)}{27\sqrt{3}}\zwe\notag\\&&+\gamma\Biggl\{\frac{32\epsilon}{3}+40A_2-\frac{382}{9}\zae+\frac{(511+53\gamma)}{6\sqrt{3}}nW_1
-\frac{(2519-24\gamma)}{27\sqrt{3}}\zwe\Biggr\}\biggr]\end{eqnarray}
\begin{eqnarray}
F_3'&=&\frac{-9}{8}\biggl[\frac{8}{3}\epsilon+\frac{203A_2}{6}-\frac{721}{54}\zae-\frac{(105+15\gamma)}{18\sqrt{3}}nW_1-\frac{(319-114\gamma)}{81\sqrt{3}}\zwe\notag\\&&+\gamma\Biggl\{2-\frac{4\epsilon}{9}-\frac{173A_2}{6}-\frac{781}{9}\zae+\frac{(197+23\gamma)}{18\sqrt{3}}nW_1
\notag\\
&&-\frac{(265-32\gamma)}{81\sqrt{3}}\zwe\Biggr\}\biggr]\end{eqnarray}
\begin{eqnarray}
F_4'&=&\frac{-3\sqrt{3}}{16}\biggl[392-\frac{532\epsilon}{3}+\frac{1918A_2}{3}-\frac{28582A_2}{9}\zae+\frac{(203+1211\gamma)}{9\sqrt{3}}nW_1\notag\\
&&+\frac{(949+4378\gamma)}{27\sqrt{3}}\zwe+28\gamma\Biggl\{\frac{108\epsilon}{7}+\frac{4037A_2}{84}-\frac{611}{21}\zae+\frac{(8397+919\gamma)}{84\sqrt{3}}nW_1
\notag\\
&&-\frac{(92266-1869\gamma)}{27\sqrt{3}}\zwe\Biggr\}\biggr]
\end{eqnarray}
\begin{eqnarray}
 F_1''&=&\frac{\zwe}{6}
\end{eqnarray}
\begin{eqnarray} F_2''&&=\frac{-9}{32}\biggl[\frac{8}{3}\epsilon+\frac{203A_2}{6}-\frac{625}{54}\zae-\frac{(105+15\gamma)}{18\sqrt{3}}nW_1-\frac{(307-114\gamma)}{81\sqrt{3}}\zwe\notag\\&&+\gamma\Biggl\{2-\frac{4\epsilon}{9}+\frac{55A_2}{2}-\frac{797}{54}\zae+\frac{(197+23\gamma)}{18\sqrt{3}}nW_1
-\frac{(211-32\gamma)}{81\sqrt{3}}\zwe\Biggr\}\biggr]\end{eqnarray}
\begin{eqnarray}
F_3''&=&\frac{-9\sqrt{3}}{16}\biggl[2-\frac{8}{3}\epsilon+\frac{55A_2}{6}-\frac{134}{3}\zae-\frac{(37+\gamma)}{18\sqrt{3}}nW_1-\frac{(93+226\gamma)}{81\sqrt{3}}\zwe
\notag\\ &&+\gamma\Biggl\{4\epsilon+\frac{169
}{27}\zae+\frac{(241+45\gamma)}{18\sqrt{3}}nW_1
-\frac{(1558-126\gamma)}{81\sqrt{3}}\zwe\Biggr\}\biggr]\end{eqnarray}
\begin{eqnarray}
F_4''&=&\frac{9}{256}\biggl[\frac{212}{3}\epsilon+\frac{2950A_2}{3}-\frac{1370A_2}{27}\zae-\frac{(771+237\gamma)}{9\sqrt{3}}nW_1-\frac{2(1907-984\gamma)}{81\sqrt{3}}\zwe\notag\\&&+28\gamma\Biggl\{\frac{11}{7}+\frac{4\epsilon}{9}-\frac{152A_2}{7}-\frac{36965}{504}\zae+\frac{(2569+277\gamma)}{252\sqrt{3}}nW_1
\notag\\&&+\frac{(22603+4396\gamma)}{1134\sqrt{3}}\zwe\Biggr\}\biggr]
\end{eqnarray}
\begin{eqnarray}
 G_1&=&\frac{-\zwe}{6}
\end{eqnarray}
\begin{eqnarray} G_2&=&\frac{3}{32}\biggl[14-\frac{16}{3}\epsilon+A_2-\frac{1367}{18}\zae+\frac{115(1+\gamma)}{18\sqrt{3}}nW_1-\frac{(863-136\gamma)}{27\sqrt{3}}\zwe\notag\\&&+\gamma\Biggl\{\frac{32\epsilon}{3}+40A_2-\frac{382 }{9}\zae+\frac{(511+53\gamma)}{6\sqrt{3}}nW_1-\frac{(2519-24\gamma)}{27\sqrt{3}}\zwe\Biggr\}\biggr]\end{eqnarray}
\begin{eqnarray}
G_3&=&\frac{3\sqrt{3}}{16}\biggl[\frac{16}{3}\epsilon+6A_2-\frac{907A_2}{18}\zae+\frac{(143+9\gamma)}{6\sqrt{3}}nW_1+\frac{(477+403\gamma)}{27\sqrt{3}}\zwe\notag\\&&+\gamma\Biggl\{14+\frac{4\epsilon}{3}+\frac{71A_2}{2}-\frac{1489}{18}\zae-\frac{(215+29\gamma)}{6\sqrt{3}}nW_1
\notag\\
&&-\frac{2(1174+169\gamma)}{27\sqrt{3}}\zwe\Biggr\}\biggr]\end{eqnarray}
\begin{eqnarray}
G_4&=&\frac{3\sqrt{3}}{256}\biggl[84+52\epsilon+212A_2-267\zae+\frac{2(299+61\gamma)}{3\sqrt{3}}nW_1-\frac{(14854+225\gamma)}{27\sqrt{3}}\zwe\notag\\&&+\gamma\Biggl\{32\epsilon+156A_2+649\zae-\frac{(562+8\gamma)}{3\sqrt{3}}nW_1+\frac{(13285+5169\gamma)}{27\sqrt{3}}\zwe\Biggr\}\biggr]
 \end{eqnarray}
\begin{eqnarray}
 G_1'&=&\frac{-\zwe}{\sqrt{3}}
\end{eqnarray}
\begin{eqnarray} G_2'&=&\frac{9}{16}\biggl[\frac{8}{3}\epsilon+\frac{203A_2}{6}-\frac{625}{54}\zae-\frac{(105+15\gamma)}{18\sqrt{3}}nW_1-\frac{(307-114\gamma)}{81\sqrt{3}}\zwe\notag\\&&-\gamma\Biggl\{2-\frac{4\epsilon}{9}-\frac{55A_2}{2}-\frac{797}{54}\zae+\frac{(197+23\gamma)}{18\sqrt{3}}nW_1
\notag\\
&&-\frac{(211-32\gamma)}{81\sqrt{3}}\zwe\Biggr\}\biggr]\end{eqnarray}
\begin{eqnarray}
G_3'&=&\frac{3\sqrt{3}}{8}\biggl[14-\frac{16}{3}\epsilon+\frac{65A_2}{6}-\frac{1439}{18}\zae+\frac{115(1+\gamma)}{18\sqrt{3}}nW_1-\frac{(941-118\gamma)}{27\sqrt{3}}\zwe\notag\\&&+\gamma\Biggl\{\frac{32\epsilon}{3}-40A_2-\frac{310}{9}\zae+\frac{(511+53\gamma)}{6\sqrt{3}}nW_1
-\frac{(251-24\gamma)}{27\sqrt{3}}\zwe\Biggr\}\biggr]\end{eqnarray}
\begin{eqnarray}
G_4'&=&\frac{-9}{128}\biggl[12\epsilon-287A_2+\frac{847A_2}{9}\zae-\frac{2(28+\gamma)}{\sqrt{3}}nW_1-\frac{4(2210-69\gamma)}{27\sqrt{3}}\zwe\notag\\&&-\gamma\Biggl\{96+\frac{152\epsilon}{3}+135A_2-\frac{2320}{9}\zae+\frac{(497-123\gamma)}{3\sqrt{3}}nW_1
\notag\\ &&-\frac{4(17697+32\gamma)}{27\sqrt{3}}\zwe\Biggr\}\biggr]
\end{eqnarray}
\begin{eqnarray}
 G_1''&=&\frac{-\zwe}{6}
\end{eqnarray}
\begin{eqnarray} G_2''&=&\frac{9\sqrt{3}}{32}\biggl[2-\frac{8}{3}\epsilon+\frac{23A_2}{3}-44\zae-\frac{(37+\gamma)}{18\sqrt{3}}nW_1-\frac{(123+349\gamma)}{3\sqrt{3}}\zwe\notag\\&&+\gamma\Biggl\{4\epsilon+\frac{88A_2}{27}+\frac{(421+45\gamma)}{18\sqrt{3}}nW_1
-\frac{(1558-126\gamma)}{81\sqrt{3}}\zwe\Biggr\}\biggr]\end{eqnarray}
\begin{eqnarray}
G_3''&=&\frac{-9}{16}\biggl[\frac{8}{9}\epsilon+\frac{203A_2}{6}-\frac{589}{54}\zae-\frac{5(51+2\gamma)}{18\sqrt{3}}nW_1-\frac{(349-282\gamma)}{81\sqrt{3}}\zwe\notag\\&&+\gamma\Biggl\{2-\frac{4\epsilon}{9}-26A_2-\frac{412}{27}\zae+\frac{(197+23\gamma)}{18\sqrt{3}}nW_1
\notag\\
&&-\frac{(211-32\gamma)}{81\sqrt{3}}\zwe\Biggr\}\biggr]\end{eqnarray}
\begin{eqnarray}
G_4''&=&\frac{-9\sqrt{3}}{256}\biggl[12+\frac{20}{3}\epsilon+76A_2-\frac{350A_2}{3}\zae+\frac{(32\gamma)}{3\sqrt{3}}nW_1-\frac{2(1529+450\gamma)}{27\sqrt{3}}\zwe\notag\\&&+\gamma\Biggl\{8\epsilon-\frac{749A_2}{3}+\frac{808}{9}\zae-\frac{(109-40\gamma)}{3\sqrt{3}}nW_1
+\frac{(35-1269\gamma)}{27\sqrt{3}}\zwe\Biggr\}\biggr]
\label{eq:g4dd}\end{eqnarray} Using transformation
$x=B_1^{\az}+B_2^{\az}$ and $y=B_1^{\za}+B_2^{\za}$ the third order
part $H_3=-L_3$ of the Hamiltonian in $I_1^{1/2},I_2^{1/2}$ is  of
the form
\begin{equation} H_3=A_{3,0}I_1^{3/2}+A_{2,1}I_1I_2^{1/2}+A_{1,2}I_1^{1/2}I_2+A_{0,3}I_2^{3/2}
\label{eq:H3}\end{equation} We can verify that in Eq.(~\ref{eq:H3}),
$A_{3,0}$ vanishes independently as in Deprit and Deprit
Bartholom\'{e}(1967). Similarly the other coefficients
$A_{2,1},A_{1,2},A_{0,3}$ are also found to be zero independently.
Hence the third order part $H_3$ of the Hamiltonian in
$I_1^{1/2},I_2^{1/2}$ is zero.

\section{Conclusion}
Using Whittaker(1965) method we have found  that the second order
part
 $H_2$ of the Hamiltonian is transformed into the normal
form $H_2=\omega_1I_1-\omega_2I_2.$ The third order part $H_3$ of the Hamiltonian in $I_1^{1/2},I_2^{1/2}$ is zero.\\

\section*{Acknowledgements} We are thankful to D.S.T. Government of India, New Delhi for sanctioning a project DST/MS/140/2K dated 02/01/2004 on this topic. We are also thankful to IUCAA Pune for providing  financial assistance for visiting library and  computer
facility.

\end{document}